\documentclass[11pt]{article}
\usepackage{latexsym,amsthm,verbatim,ifthen,graphicx,color,mathrsfs}
\usepackage[all]{xy}
\usepackage{color}

\usepackage{color}

\usepackage{lscape}
 \usepackage[ansinew]{inputenc}
 \usepackage{multicol}
 \usepackage[all]{xy}\xyoption{poly}
 \usepackage{mathrsfs}
\usepackage[psamsfonts]{eucal}
\usepackage{amssymb, amsmath}
\usepackage{geometry}
\usepackage{enumerate}
\usepackage{float}

\newtheorem{theorem}{Theorem}[section]

\newtheorem{corollary}[theorem]{Corollary}
 \newtheorem{lemma}[theorem]{Lemma}
 \newtheorem{proposition}[theorem]{Proposition}
 \theoremstyle{definition}

\def\quic{{\mathscr C}}

\newcommand{\bz}{\mathbb Z}

\def\bq{\mathbb{Q}}

\def\L{\mathbb{L}}

\def\Q{\mathbb{Q}}

\def\hL{{\widehat{\mathbb L}}}

    \newcommand{\lasu}{{\mathfrak{L}}}

\newcommand{\lpol}{\operatorname{\text{\L}}}

 \newcommand{\lib }{\mathbb{L}}

  \newcommand{\Hom}{\operatorname{\text{\rm Hom}}}
    
\newcommand{\catcdga}{\operatorname{{\bf cdga}}}

\newcommand{\catrdgl}{\operatorname{{\bf dgl_1}}}

\newcommand{\catcdgl}{\operatorname{{\bf cdgl}}}
\newcommand{\catedgl}{\operatorname{{\bf edgl}}}

\newcommand{\catsgp}{\operatorname{{\bf sgp_0}}}
\newcommand{\catshopf}{\operatorname{{\bf sch_0}}}
\newcommand{\catss}{\operatorname{{\bf sset}}}

\newcommand{\catssr}{\operatorname{{\bf sset_{1} }}}
\newcommand{\catssrq}{\operatorname{{\bf sset_1^\bq}}}
\newcommand{\catslie}{\operatorname{{\bf sla_1}}}

\newcommand{\Ho}{\operatorname{{\rm Ho}}}

    \newcommand{\id}{\operatorname{{\rm id}}}
        \newcommand{\ad}{\operatorname{{\rm ad}}}

 \newcommand{\MC}{\operatorname{{\rm MC}}}

\newcommand{\mc}{{\MC}}

  \newcommand{\otimesc}{\widehat{\otimes}}

\newcommand{\hocolim}{\underleftarrow
{\operatorname{\mathstrut hocolim}}}

   \newcommand{\libc}{{\widehat\lib}}

\usepackage{hyperref}
\usepackage{physics}
\hypersetup{
    colorlinks=true,
    linkcolor=blue,
    filecolor=magenta,
    urlcolor=cyan,
    citecolor=blue,
}

\newcommand{\tensor}{\otimes}

\newtheorem{example}[theorem]{Example}				
\theoremstyle{remark}
\newtheorem{remark}[theorem]{Remark}

\usepackage{stmaryrd}

\renewcommand{\P}{\CMcal{P}}
\usepackage{ textcomp }
\newcommand{\cdgl}{\operatorname{\bf cdgl}}
\newcommand{\dbar}{\bar d}
\newcommand{\sbar}{\bar s}
\renewcommand{\[}{\llbracket}
\renewcommand{\]}{\rrbracket}

 \newcommand{\asu}{{\mathscr{A}}}

\begin{document}

\title{All known realizations of complete Lie algebras coincide}

\author{Yves F\'elix, Mario Fuentes and Aniceto Murillo\footnote{The  authors have been partially supported by the MICINN grant PID2020-118753GB-I00 of the Spanish Government.\hfill\break MSC2020: 55P62, 17B55.}}

\maketitle

\begin{abstract}
We prove that for any reduced differential graded Lie algebra  $L$, the classical Quillen geometrical realization $\langle L\rangle_Q$ is homotopy equivalent to the realization $\langle L\rangle=\Hom_{\cdgl}(\lasu_\bullet,L)$ constructed via the cosimplicial free complete differential graded Lie algebra $\lasu_\bullet$. As the latter is  a deformation retract of the Deligne-Getzler-Hinich realization $\mc_\bullet(L)$   we deduce that, up to homotopy, all known topological realization functors of complete differential graded Lie algebras coincide. Immediate consequences of our main result include an elementary proof of the Baues-Lemaire conjecture and the description of the Quillen realization as a representable functor.
\end{abstract}

\section*{Introduction}
In  \cite{qui}, Quillen constructed a geometrical realization functor
$$
\langle\,\cdot\,\rangle_Q\colon \catrdgl\longrightarrow\catssr
$$
from the category of simply connected differential graded Lie algebras to the category of reduced simplicial sets.  This was the starting point of rational homotopy theory, from the Lie approach, as this functor induces an equivalence between the corresponding homotopy categories when considering rational reduced simplicial sets. Later on \cite{getz,hi}, the {\em Deligne-Getzler-Hinich groupoid} functor
$$
  \mc_\bullet\colon\catcdgl\longrightarrow \catss
  $$
  was defined in the category of complete differential graded Lie algebras (cdgl's henceforth). Given $L$ a cdgl, $ \mc_\bullet(L)=\mc(\asu_\bullet\otimesc L)$ is the simplicial set of Maurer-Cartan elements of the simplicial cdgl $\asu_\bullet\otimesc L$ in which $\asu_\bullet$ denotes the simplicial commutative differential graded algebra  of PL-differential forms on the standard simplices.

 Finally, see \cite{bfmt0}, there is a realization  functor for cdgl's based on a quite geometrical cosimplicial cdgl $\lasu_{\bullet}$,
 $$
 \langle\,\cdot\,\rangle\colon \catcdgl\longrightarrow \catss,\quad \langle L\rangle=\Hom_{\cdgl}(\lasu_\bullet,L).
 $$
In  \cite[Theorem 0.2]{ter} and \cite[Theorem  5.2]{ni1} it is shown that, for any cdgl $L$,
 $\langle L\rangle$ is simplicially isomorphic to $\gamma_\bullet (L)$, the {\em nerve of $L$} \cite[\S5]{getz},  which is  a deformation retract of  $\mc_\bullet(L)$.

  Here, we close up the circle and prove:

   \begin{theorem}\label{intro} For any simply connected differential graded Lie algebra $L$,
   $$
   \langle L\rangle\simeq\langle L\rangle_Q.
   $$
   \end{theorem}
    The first immediate consequence  is that the functors induced in the respective homotopy categories by the global model and realization functors, see \cite{bfmt0},
$$
\xymatrix{ \catss& \catcdgl \ar@<1ex>[l]^(.50){\langle\,\cdot\,\rangle}
\ar@<1ex>[l];[]^(.50){\lasu}\\}
$$
extend (in a unique way) the classical Quillen equivalence between the homotopy categories of rational reduced simplicial sets and that of simply connected differential graded Lie algebras. In particular, and under no restriction, the
    Quillen original realization functor  is representable by the cosimplicial Lie algebra $\lasu_\bullet$ and thus, it can be finally regarded as the Eckmann-Hilton dual of the  realization functor of cdga's, co-representable by $\asu_\bullet$. That is a question which has puzzled rational homotopists since the birth of both, the Sullivan and Quillen approaches to rational homotopy theory.

Theorem \ref{intro} was already known for simply connected dgl's of finite type, see \cite[Corollary 11.17]{bfmt0} but its proof heavily relied in the Baues-Lemaire Conjecture \cite[Conjecture 3.5]{baule} proved in \cite{ma}. The opposite procedure is now available. Namely, the self-contained proof of Theorem \ref{intro} under no restriction let us trivially reprove this conjecture  by which, given $L$ the minimal Quillen model of a simply connected complex of finite type $X$, the cdga $\mathscr C^*(L)$ of cochains on $L$ has the homotopy type of  the Sullivan minimal model of $X$.

For completeness, we also collect extensions of this conjecture to the non finite type and/or non simply connected case and describe, in  this extended scenario, the result of  composing  the Sullivan realization functor of commutative differential graded algebras (cdga's henceforth) \cite{su},
$$
  \langle\,\cdot\,\rangle_S \colon \catcdga\longrightarrow \catss,\quad   \langle A\rangle_S=\Hom_{\catcdga}(A,\asu_\bullet),
   $$
   and the appropriate extension  of the Chevalley-Eilenberg cochain functor adapted to the specific class of considered differential graded Lie algebras.

\section{Preliminaries}\label{prelimi}

 Unless explicitly stated otherwise, any algebraic object  is assumed to be $\bz$-graded and having $\bq$ as coefficient field. We denote in bold the category containing such object.

 A \emph{differential graded Lie algebra} (dgl henceforth)  is a graded vector space $L$ endowed with a {\em Lie bracket}  $[\,\,,\,]$ satisfying graded antisymmetry and Jacobi identity, and a linear derivation   $d$  of degree $-1$ such that $d^2= 0$. An element $a\in L_{-1}$ is {\em Maurer-Cartan} (MC element hereafter) if $da=-\frac{1}{2}[a,a] $. Given an MC element $a\in L_{-1}$ the map $d_a=d+\ad_a$ is a new differential on $L$.  The {\em component} of $L$ at $a$ is the connected sub dgl $L^a$ of $(L,d_a)$ given by
$$
L^a_p=\begin{cases} \ker d_a&\text{if $p=0$},\\ \,\,\,L_p&\text{if $p>0$}.\end{cases}
$$
 A dgl $L$ is called \emph{free} if it is free as a Lie algebra, that is, $L = \mathbb L(V)$ for some graded vector space $V$.  A dgl $L$ is {\em simply connected or reduced} if it is concentrated in positive degrees, $L=\oplus_{p\ge1}L_p$.

 A {\em complete differential graded Lie algebra}, cdgl henceforth,  is a dgl $L$ equipped with a  decreasing filtration of differential ideals,
 $$L=F^1\supset\dots \supset F^n\supset F^{n+1}\supset\dots,$$
 with $[F^p,F^q]\subset F^{p+q}$ for $p,q\geq 1$ and such that the natural map
$$
L\stackrel{\cong}{\longrightarrow}\varprojlim_n L/F^n
$$
is a dgl isomorphism. The lower central series of a given dgl,
$$
L^1\supset\dots\supset L^n\supset L^{n+1}\supset\dots,\quad L^{1}= L,\quad L^{n}= [L, L^{n-1}],
$$
is always a filtration. A {\em morphism} $f\colon L\to L'$ between cdgl's is a filtration preserving dgl morphism. Note that any simply connected dgl is always complete.

Given  a free Lie algebra  $\lib(V)$ we write
 $$
\libc(V)=\varprojlim_n\lib(V)/\lib(V)^n$$
which is complete with respect to the filtration
 $ F^n=  \ker(\libc(V)\to \lib(V)/\lib(V)^n)$, $n\ge 1$.

The homotopy theory of cdgl's, for which we refer to \cite{bfmt0}, is based on the pair of
 of adjoint functors, {\em (global) model} and {\em realization},
\begin{equation}\label{pair}
\xymatrix{ \catss& \catcdgl \ar@<0.75ex>[l]^(.50){\langle\,\cdot\,\rangle}
\ar@<0.75ex>[l];[]^(.50){\lasu}\\},
\end{equation}
which highly rely on the cosimplicial cdgl
$\lasu_\bullet=\{\lasu_n\}_{n\ge 0}
$, see \cite[Chapter 6]{bfmt0}.
 For each $n\ge 0$, ,
$$
\lasu_n=\bigl(\libc(s^{-1}\Delta^{n}),d)
$$
where
$s^{-1}\Delta^{n}$  denotes the desuspension  of the non-degenerate simplicial chains  on the simplicial set $\underline\Delta^n$. That is, for any $p\ge 0$, a generator of  degree $p-1$ of $s^{-1}\Delta^n$ can be written as $a_{i_0\dots i_p}$ with $0\le i_0<\dots<i_p\le n$. The cofaces and codegeneracies in $\lasu_\bullet$ are induced by those on the cosimplicial chain complex $s^{-1}\Delta^{n}$:

For each $0\le i\le n$, the  map $\delta_i   \colon \{0,\dots, n-1\}\to \{0, \dots, n\}$,
$$
\delta_i(j)=\begin{cases}\,\,j\,\,&\text{if $j<i$,}\\ j+1\,\,&\text{if $j\ge i$},
\end{cases}
$$
defines the cdgl morphism,
$$
\delta_i\colon \libc(s^{-1}\Delta^{n-1})\longrightarrow \libc(s^{-1}\Delta^n),\quad
\delta_i(a_{\ell_0\dots \ell_p}) = a_{\delta_i(\ell_0)\dots \delta_i(\ell_p)}.
$$
On the other hand, the map $\sigma_i   \colon \{0,\dots, n+1\}\to \{0, \dots, n\}$,
$$\sigma_i(j) = \begin{cases} \,\, j&\text{if $j\le i$},\\ j-1&\text{if $j> i$},
\end{cases}
$$
also defines the cdgl morphism $\sigma_i\colon \libc(s^{-1}\Delta^{n+1})\to  \libc(s^{-1}\Delta^{n})$,
$$
\sigma_i (a_{\ell_0\dots \ell_q}) = \left\{\begin{array}{ll}
a_{\sigma_i (\ell_0)\dots  \sigma_i (\ell_q)}   & \mbox{if } \sigma_i (\ell_0)<\dots <\sigma_i (\ell_q),
\\
\qquad 0 & \mbox{otherwise.}
\end{array}\right.
$$
The differential $d$ on each $\lasu_n$ is the only one (up to cdgl isomorphism) satisfying:

 \begin{itemize}
 \item[(1)] For each $i=0,\dots,n$, the generators $a_0,\dots,a_n\in s^{-1}\Delta^{n}$, corresponding to vertices, are MC elements.

     \item[(2)] The linear part of $d$ is induced by the boundary operator of $s^{-1}\Delta^{n}$.
     \item[(3)] The cofaces and codegeneracies are cdgl morphisms.
     \end{itemize}
The realization of any cdgl $L$ is defined as the simplicial set,
$$
\langle L\rangle={\Hom_{\cdgl}}(\lasu_\bullet,L),\quad d_i=\Hom_{\cdgl}(\delta_i,L),\quad s_i=\Hom_{\cdgl}(\sigma_i,L).
$$
If $L$ is connected, that is, non negatively graded, then $\langle L\rangle$ is a connected simplicial set and for any $n\ge 1$, the map
\begin{equation}\label{homotopy}
     \rho_n\colon \pi_n\langle L\rangle\stackrel{\cong}{\longrightarrow} H_{n-1}(L),\quad \rho_n[\varphi]=[\varphi(a_{0\dots n})],
\end{equation}
    is a group isomorphism \cite[Theorem 7.18]{bfmt0}. Here, the group law in $H_0(L)$ is given by the Baker-Campbell-Hausdorff  product (BCH product henceforth).

On the other hand, in the seminal paper \cite{qui}, for which we refer for details of what follows, D. Quillen introduced a couple of  functors,
$$
\xymatrix{  \catssr \ar@<0.75ex>[r]^-{\lambda} &\catrdgl \ar@<0.75ex>[l]^(.40){\langle\,\cdot\,\rangle_Q} }
$$
which are  the composition of the following pairs of adjoint functors (the upper arrow denotes left adjoint),
$$
\xymatrix{  \lambda\colon \catssr \ar@<0.75ex>[r]^-{G} &\catsgp \ar@<0.75ex>[l]^(.40){W}
\ar@<0.75ex>[r]^-{\widehat\bq} &\catshopf \ar@<0.75ex>[l]^(.46){\mathcal G}
\ar@{<-}@<0.75ex>[r]^-{\widehat U}  &\catslie \ar@{<-}@<0.75ex>[l]^(.46){\mathcal P}
\ar@{<-}@<0.75ex>[r]^-{N^*}  &\catrdgl \ar@{<-}@<0.75ex>[l]^(.65){N}\colon \langle \,\cdot\, \rangle_Q.
}
$$
Here,  $\catssr $, $\catsgp$, $\catshopf$ and $\catslie$ denote respectively the categories of {\em reduced simplicial sets} (those with only one simplex in dimensions $0$ and $1$), {\em connected simplicial groups} ($G_0=\{1\}$), {\em connected complete Hopf algebras} ($A_{<0}=0$ and $A_0=\bq$), and {\em reduced simplicial Lie algebras}. Each of these  pairs induces
Quillen equivalences on  the corresponding homotopy categories when localizing on  rational weak homotopy equivalences in $\catssr$ and $\catsgp$,  on weak homotopy equivalences in $\catshopf$ and $\catslie$, and on quasi-isomorphisms in $\catrdgl$ \cite[Thm. I]{qui}.

 Recall that, given $X$ a reduced simplicial set,  $G(X)$ is the Kan simplicial group where, for any $n\ge 0$, $(GX)_n$ is the free group generated by $X_{n+1}\setminus s_0X_n$ or equivalently, the quotient group
$(GX)_n=\text{Free}(X_{n+1})/(s_0X_n)$. Faces and degeneracies are defined as:
$$\begin{array}{lll}
\sbar_i\colon  (GX)_n\to (GX)_{n+1}, & \sbar_ix=s_{i+1}x, \quad x\in X_{n+1},\quad i=0,\dots, n,\\
 \\
\dbar_i\colon (GX)_n\to (GX)_{n-1}, & \dbar_ix=d_{i+1}x, \quad x\in X_{n+1},\quad i=1,\dots, n,\\
\\
\dbar_0\colon (GX)_n\to (GX)_{n-1}, & \dbar_0x= (d_0x)^{-1}(d_{1}x),\quad x\in X_{n+1}.\\
\end{array}$$
Next, see for instance \cite[Proposition 8.4.1]{fres}, the composition $\P\hat\Q GX$ is the simplicial free complete Lie algebra
$$\lpol_n=\hL(X_{n+1}/s_0X_n),\quad n\ge0,$$
in which $X_n$ denotes here the vector space generated by the $n$-simplices of $X$. We denote by $\[\cdot,\cdot\]$ the Lie bracket in this simplicial Lie algebra.
The faces and degeneracies $\dbar_i$'s  and $\sbar_i$'s acts as above on generators and are extended as Lie morphisms. Note that, since the multiplication on $GX$ is taking to the BCH product on $\lpol$,
$$\dbar_0 x=(-d_0x)\ast (d_{1}x),\quad x\in X_n,$$
where $\ast$ denotes the BCH product.

Finally $\lambda(X)= N\lpol$ is the reduced dgl given by the {\em normalized chain complex} on $\lpol$,
$$(N\lpol)_n=\bigcap_{i=1}^n \ker(\dbar_i\colon\lpol_n\to \lpol_{n-1}),\quad n\ge 1,\quad d=\dbar_0,$$
endowed with the bracket induced by the simplicial  Eilenberg-Zilber formula: given $\omega\in (N\lpol)_n$ and $\omega'\in (N\lpol)_m$,
$$[ \omega,\omega']=\sum_{(\mu,\nu)\in S_{n,m}^0} \varepsilon_{\mu,\nu} \[\sbar_{\nu_m} \sbar_{\nu_{m-1}} \cdots \sbar_{\nu_1} \omega, \sbar_{\mu_n} \sbar_{\mu_{n-1}}\cdots \sbar_{\mu_1}\omega'\],$$
where $ S_{n,m}^0$ denotes the $(n,m)$-shuffles of $\{ 0,1,\dots, n+m-1\}$  and $\varepsilon_{\mu,\nu}$ is the sign of such shuffle.

In general, given a simplicial vector space $V$ we denote by $C(V)$ its simplicial  chain complex in which $C_n(V)=V_n$ and the boundary operator is given by $\sum_i(-1)^id_i$. Given another simplicial vector space $W$, there are natural chain maps, Eilenberg-Zilber and Alexander-Whitney (see for instance \cite[\S29]{may}),
$$
\xymatrix{C(V)\otimes C(W)& C(V\otimes W), \ar@<0.75ex>[l]^(.44){\Delta}
\ar@<0.5ex>[l];[]^(.54){\nabla}\\}
$$
in which $V\otimes W$ is the simplicial tensor product, given by
$$
\nabla(\alpha\tensor \beta)= \sum_{ (\mu,\nu)\in  S_{n,m}^0 } \varepsilon_{\mu,\nu}\,
s_{\nu}\alpha \, \tensor \, s_{\mu}\beta,\qquad \alpha\in V_n,\quad \beta\in W_m,$$
where $s_\nu$ and $s_\mu$ stand for
$s_{\nu_m}s_{\nu_{m-1}}\cdots s_{\nu_{1}}$ and $s_{\mu_n}s_{\mu_{m-1}}\cdots s_{\mu_{1}}$, and
 $$
\Delta(\alpha\tensor\beta)=
\sum_{k=0}^n d_{k+1}^{n-k} \alpha \; \tensor \; d_0^{k} \beta=
\sum_{k=0}^n d_{k+1}d_{k+2}\cdots d_{n} \alpha\;\tensor\; d_0d_1\cdots d_{k-1} \beta,\quad\alpha\otimes\beta\in V_n\otimes W_n.
$$
These maps restrict to the corresponding normalized chain complexes
$$
\xymatrix{NV\otimes NW& N(V\otimes W) \ar@<0.75ex>[l]^(.47){\Delta}
\ar@<0.5ex>[l];[]^(.50){\nabla}\\}
$$
where they satisfy
\begin{equation}\label{retracto}
\Delta\nabla=\id_{NV\otimes NW},\qquad \nabla\Delta\simeq\id_{N(V\otimes W)}.
\end{equation}
With this notation,  the Lie bracket  on $\lambda(L)=N\lpol$ is given by
\begin{equation}\label{liebracket}
[ \cdot, \cdot]=\[\cdot,\cdot\]\circ\nabla\colon  N\lpol \tensor N\lpol\to N\lpol.
\end{equation}

\section{The proof and consequences}

In view of (\ref{homotopy}), if  $L$ is a simply connected dgl then $\langle L\rangle$ is a rational simplicial set. As the Quillen functors induce equivalences between the homotopy categories of reduced rational simplicial sets and that of simply connected dgl's, Theorem \ref{intro} is an immediate consequence of the following, more general, result:

\begin{theorem}\label{main}
Given a simply connected dgl $L$ there is a surjective quasi-isomorphism,
$$
\Phi\colon\lambda\langle L\rangle\stackrel{\simeq}{\longrightarrow} L.
$$
\end{theorem}

\begin{proof}
Write $X_\bullet=\langle L\rangle=\Hom_{\catcdgl}(\lasu_\bullet,L)$. We first define a linear map
$$
\Phi\colon C(\lpol)\longrightarrow L
$$
recursively on each $\lpol_n^m=\libc^m(X_{n+1}/s_0X_n)$,    being $n\ge 1$ the simplicial degree and $m\ge1$ the word length:

If $\varphi\colon\lasu_2\to L\in X_2$ is an indecomposable of $\lpol_1$ define $\Phi(\varphi)=-\varphi(a_{012})\in L_1$. Observe that if $\varphi=s_0\eta$, then $\Phi(\varphi)=-\eta\bigl(\sigma_0(a_{012})\bigr)=-\eta(0)=0$. Set $\Phi(\lpol_1^{\ge2})=0$.

Assume $\Phi$ built on $\lpol_{<n}$ and define
$$
\Phi\colon \lpol_n\longrightarrow L
$$
as follows: again, for an indecomposable $\varphi\colon\lasu_{n+1}\to L\in X_{n+1}$  define $$\Phi(\varphi)=(-1)^n\varphi(a_{0\dots n+1}).$$
Once more, if $\varphi=s_0\eta$, then $\Phi(\varphi)=0$ and thus $\Phi$ is well defined on $X_{n+1}/s_0X_n$. Now, if $\Phi$ has been defined for $\alpha,\beta\in \lpol_n$, set
\begin{equation}\label{eq:Phi}
\Phi\[\alpha,\beta\]=[\cdot,\cdot]\circ(\Phi\tensor\Phi)\circ\Delta(\alpha\tensor\beta).
\end{equation}
As $\Phi$ increases the bracket length on each $\lpol_n$, it is well defined as terms of high enough word length of any series in $\lpol_n$ are sent to $0$.

Finally, recall that $\lambda(X)=N\lpol$ and  define
$$
\Phi\colon\lambda\langle L\rangle\stackrel{}{\longrightarrow} L
$$
as the restriction of $\Phi$ to $N\lpol\subset C(\lpol)$.

\medskip

{\em (i)} $\Phi\colon N{\lpol}\to L$ {\em is a morphism of Lie algebras:} Indeed, by (\ref{liebracket}), (\ref{eq:Phi}) and (\ref{retracto}) respectively, we have the following identities:
$$\Phi\circ[ \cdot,\cdot]=\Phi\circ \[\cdot,\cdot\]\circ \nabla= [\cdot,\cdot]\circ (\Phi\tensor \Phi)\circ \Delta\circ\nabla=[\cdot,\cdot] \circ (\Phi\tensor\Phi).$$

{\em (ii)} $\Phi\colon N{\lpol}\to L$ {\em is a surjective chain map:} We first check that it is enough to show that $\Phi d=d\Phi$ for the indecomposable elements  of $N\lpol$ (with respect to the Lie bracket $\[\cdot,\cdot\]$). Assume that this is the case and suppose, inductively, that  $\Phi(d \alpha)=d \Phi(\alpha)$ and $\Phi(d \beta)=d\Phi(\beta)$ for some elements $\alpha,\beta\in (N\lpol)_n$. Then,
$$
\begin{aligned}
\Phi(d \[\alpha,\beta\])&=\Phi\circ \[\cdot,\cdot\] \circ (d \tensor d)  (\alpha\tensor \beta)
=
[\cdot,\cdot] \circ (\Phi\tensor\Phi)\circ \Delta \circ (d \tensor d)  (\alpha\tensor \beta)\\
&=
[\cdot,\cdot]\circ (\Phi\tensor \Phi)\circ (d\tensor \id+ (-1)^n \id\tensor d)\circ \Delta(\alpha\tensor \beta)\\
&=
[\cdot,\cdot]\circ (d \tensor \id+(-1)^n\id\tensor d)\circ (\Phi\tensor \Phi)\Delta(\alpha\tensor \beta)\\
&=d\circ[\cdot,\cdot]\circ (\Phi\tensor\Phi)\Delta(\alpha\tensor \beta)=d\Phi\[\alpha,\beta\],
\end{aligned}
$$
where, for the above identities we have used, respectively,  the following facts: $d=\dbar_0$ in $N\lpol$ so it commutes with $\[\cdot,\cdot\]$; formula (\ref{eq:Phi}); $\Delta$ is a chain map; induction hypothesis;  $d$ is a derivation on $L$; and again  formula (\ref{eq:Phi}).

We next prove that   $\Phi$ commutes with the differential for indecomposable elements of $N\lpol$. Note that, for a fixed degree $n\ge 1$, these are the vector subspace of $X_{n+1}$ generated by the non degenerate simplices.

Choose $\varphi\colon\lasu_{n+1}\to L$ such a generator so that $\dbar_i\varphi=0$ for $i=1,\dots n$, which is equivalent to $\varphi\circ\delta^{i}=0$ for $i=2,\dots, n+1$. This implies that  the only possible non-zero images of generators are
\begin{equation}\label{solo}
\varphi(a_{01\dots n+1})=x_{01}, \;\;\; \varphi(a_{02\dots n+1})=x_0,
\;\;\;
\varphi(a_{12\dots n+1})=x_1, \;\;\; \varphi(a_{23\dots n+1})=x_2,
\end{equation}
which are elements in $L$ of degree $n, n-1, n-1$ and $n-2$ respectively. Recall  that, in $\lasu_{n+1}$,
$$
da_{01\dots n+1}=
\sum_{i=0}^{n+1} (-1)^i a_{0\dots \hat i \dots n+1}+\Omega_n
$$
where $\Omega_n$ is a decomposable element of degree $n-1$. By degree reasons, one easily sees that $\varphi(\Omega_n)=0$ except, at most, for $n=3$ and as long as $\Omega_3$ contains $[a_{234},a_{234}]$ as a summand.  The technical Lemma \ref{lema} below shows that this is not the case for an appropriate choice of the differential in $\lasu_4$ and therefore,
$$\varphi(d a_{01\dots n+1})=\varphi(a_{12\dots n+1})-\varphi(a_{02\dots n+1}).$$
Another easy inspection let us also write
$$\varphi(d a_{02\dots n+1})=\varphi(d a_{12\dots n+1})=\varphi(a_{2\dots n+1})$$
and we conclude that
\begin{equation}\label{determina}
d x_{01}=x_1-x_0\quad\text{and}\quad d x_0=d x_1=x_2.
\end{equation}

On the one hand,
$$
d\Phi(\varphi)=(-1)^ndx_{01}=(-1)^n(x_1-x_0).
$$

On the other hand,
$$
\begin{aligned}
\Phi(d\varphi)&=\Phi(\dbar_0\varphi)=\Phi\bigl((-d_0\varphi)\ast (d_1\varphi)\bigr)=\Phi(d_1\varphi-d_0\varphi)+\Phi(\Psi)\\
&=(-1)^{n-1}(x_0-x_1)+\Phi(\Psi)=(-1)^n(x_1-x_0)+\Phi(\Psi),
\end{aligned}
$$
where $\Psi$ is an infinite sum of decomposable $\[\cdot,\cdot\]$-brackets whose letters are either $d_0\varphi$ and $d_1\varphi$. Taking into account  how $\Phi$ operates in $\[\cdot,\cdot\]$-brackets, see (\ref{eq:Phi}),  the only possible non-zero values of $\varphi$, see (\ref{solo}), and degree arguments, one checks that $\Phi(\Psi)=0$ and therefore $d\Phi(\varphi)=\Phi(d\varphi)$.

Finally, to see that $\Phi$ is surjective, given $x\in L_n$ consider the  decomposable element $\varphi$ in $(N\lpol)_n$ defined by
$$
\varphi(a_{01\dots n+1})=(-1)^n x, \quad\varphi(a_{12\dots n+1})=(-1)^n dx, \quad\varphi(a_{02\dots n+1})=\varphi(a_{23\dots n+1})=0,
$$
In view of (\ref{solo}) and (\ref{determina}) $\varphi$ is well defined and $\Phi(\varphi)=x$.

\medskip

{\em (iv)} $\Phi\colon N{\lpol}\to L$ {\em is a quasi-isomorphism:}  Recall our notation $X=\Hom_{\catcdgl}(\lasu_\bullet, L)$ and consider, for any $n\ge 1$, the composition
$$
\xymatrix{ \rho_{n+1}\colon \pi_{n+1}(X) \ar[r]^(.55){\gamma}_(.55){\cong} &H_n(N\lpol) \ar[r]^{H_n(\Phi)}&H_n(L) }
$$
where $\gamma\colon\pi_{n+1}(X)\stackrel{\cong}{\longrightarrow}H_n(\lambda X)$ is the isomorphism induced by $\lambda$ in homotopy and homology groups respectively. Recall once again that $\lambda=N\P\hat\Q G$ and through this composition  the class of an $(n+1)$-simplex $\varphi$ representing an element in $\pi_{n+1}(X)$ is taken precisely to the homology class of the indecomposable element $\varphi$ whose degree has been shifted by 1.

Indeed, $[\varphi]\in\pi_{n+1}(X)$ is trivially sent to $[\varphi]\in \pi_n(GX)$. By \cite[\S3]{qui} this is taken to the homotopy class of the generator $\varphi$ of the simplicial free Lie algebra $\lpol_n$. Finally, by (4.1) of $\cite[\S4]{qui}$ this homotopy class is sent to the homology class $[\varphi]\in (N\lpol)_n$

Hence $\rho_{n+1}[\varphi]=(-1)^n[\varphi(a_{01\dots n+1})]$ which is, up to the sign, the isomorphism  (\ref{homotopy}) and we conclude that $H_n(\Phi)$ is an isomorphism for any $n\ge 1$.

\end{proof}

\begin{lemma}\label{lema} The term $\lasu_4$ of the cosimplicial cdgl $\lasu_\bullet$ can be endowed with a differential such that the term $[a_{234},a_{234}]$ does not appear in $da_{01234}$.
\end{lemma}

\begin{proof}
Write $V=s^{-1}\Lambda^4_4\subset s^{-1}\Delta^4$ where $\Lambda^4_4$ denotes the 4th horn of $\Delta^4$.
Recall from the proof of \cite[Theorem 6.7]{bfmt0} that in $\lasu_4$ the differential of the top generator is given by
$$d_{a_0} a_{01234}=a_{0123}-\Gamma$$ where $\Gamma\in \hL(V)$ is any  solution of the equation
\begin{equation}\label{recurso}
d_{a_0}a_{0123}=\partial_{a_0}\Gamma,
\end{equation}
which always exists as $H(\libc(V),d_{a_0})=0$.
Write $\Gamma=\sum_{j\ge 1}\Gamma_j$, $d_{a_0}a_{0123}=\sum_{k\ge 1} \omega_k$, where the subscripts denotes  word lengths, and $d_{a_0}=\sum_{i\ge 1}d_i$ where each $d_i$ (not to be confused with face maps!) increases the word length by $i-1$. Then, equation(\ref{recurso}) becomes the following family of equations in $\libc(V)$:
$$
\sum_{j+i=k+1} d_i\Gamma_j=\omega_k,\qquad k\ge1,
$$
which can be recursively solved  by choosing $\Gamma_k$ so that
\begin{equation}\label{eq:Gamma}
d_1\Gamma_k=\omega_k-\partial_2\Gamma_{k-1}-\dots-\partial_{k+1}\Gamma_1.
\end{equation}
Note that such an element exists since  $H(\libc(V),d_1)=0$.

Next, we modify such a given solution as follows. With the help of a computer we  find the  element in $\mathbb{L}^2(V)$,
$$\gamma=2[a_{24}-a_{23}-a_{34},a_{0234}]
-[a_{023},a_{023}]+2[a_{023},a_{024}]-2 [ a_{023}, a_{034}] - [a_{024},a_{024}] + $$
$$
+2 [a_{024},a_{034} ] -[a_{034},a_{034}] +[a_{234},a_{234}],
$$
which contains the summand  $[a_{234},a_{234}]$ and satisfies $d_1\gamma=0$ as it can be easily checked. Hence, for a good choice of $\lambda\in \Q$, the term  $[a_{234},a_{234}]$ does not appear in the element
$\Gamma_2+\lambda \gamma$. Furthermore,
$$d_1(\Gamma_2+\lambda\gamma)=d_1\Gamma_2=\omega_2-d_2\Gamma_1$$
so  equation (\ref{eq:Gamma}) holds for $\Gamma'_1=\Gamma_1$ and $\Gamma'_2=\Gamma_2+\lambda \gamma$. Solving, as stated,  equation $(\ref{eq:Gamma})$ recursively for $k\geq 3$, we get a sequence $\Gamma'_k$ of elements in $\mathbb{L}^k(V)$ so that $\Gamma'=\sum_{k=1}^\infty \Gamma'_k$ is another solution of (\ref{recurso}). Thus, in $\lasu_4$, the differential of the top generator, which is now $d_{a_0} a_{01234}=a_{0123}-\Gamma'$, does not contain the term  $[a_{234},a_{234}]$.
\end{proof}

Recall from \cite[Chapter 8]{bfmt0} that (\ref{pair}) constitutes a Quillen pair with respect to the usual model category on $\catss$ and the one on $\catcdgl$ given in \cite[\S8.1]{bfmt0}. Then, the first consequence of Theorem \ref{intro} reads:

\begin{corollary}
The adjoint functors
$$
\xymatrix{ \Ho\catss& \Ho\catcdgl \ar@<1ex>[l]^(.50){\langle\,\cdot\,\rangle}
\ar@<1ex>[l];[]^(.50){\lasu}\\}
$$
extend the Quillen equivalences
$$
\xymatrix{ \Ho \catssrq \ar@<0.75ex>[r]^-{\lambda} &\Ho\catrdgl \ar@<0.75ex>[l]^(.40){\langle\,\cdot\,\rangle_Q}. }
$$
\end{corollary}
Here, $ \Ho \catssrq$ denotes the homotopy category of rational reduced simplicial sets.
\begin{proof} On the one hand, by Theorem \ref{intro}, $\langle L\rangle\simeq \langle L\rangle_Q$ for any simply connected dgl $L$.

On the other hand, for any reduced simplicial set $X$, consider $\lasu_X^a$ where $a$ is the only $0$-simplex of $X$ (see \cite[Chapter 6 and 8]{bfmt0}). We finish by checking that $\lasu_X^a\simeq \lambda(X)$. For it, write $X$ as the homotopy colimit $\hocolim_i X_i$ where $X_i$ are finite type subsimplicial sets of $X$. As $\lambda$ and $\lasu$ are, respectively, equivalence and left adjoint functor between the corresponding homotopy categories, they both preserve homotopy colimits. Hence,
$$
\lambda(X)=\lambda(\hocolim_i\,X_i)\simeq \hocolim_i\,\lambda(X_i)\simeq \hocolim_i\,\lasu_{X_i}^a\simeq\lasu_X^a,
$$
where the third identity is \cite[Theorem 10.2]{bfmt0}.
\end{proof}

\begin{remark}\label{unico}
In fact, the model and realization functors constitute the only Quillen pair extending, up to homotopy, the classical Quillen functors in the following sense: it is known that a Quillen model $(\lib(V),d)$ of a given reduced simplicial set $X$ can be chosen to be generated by the desuspension of non degenerate simplices  and where the differential reflects the simplicial structure of $X$. In particular, if we write $d=\sum_{i\ge 1}d_i$ where $d_i$ increases the bracket length by $i-1$, it follows that $d_1\colon V\to V$ is the desuspension of the chain differential on the non degenerate simplices and $d_2\colon V\to V\otimes V$ corresponds to  an approximation of the diagonal. If we now allow simplices in degree 0 and 1, then vertices of $X$ must correspond to MC elements of $V$. As asserted in \cite[Proposition 7.8]{bfmt0} there is a unique differential $d$ in $\libc(V)$ fulfilling this properties and $\lasu(X)=(\libc(V),d)$ is the model functor whose (unique) adjoint is necessarily the realization functor $\langle\,\cdot\,\rangle$.
\end{remark}

\begin{corollary} {\em (Baues-Lemaire conjecture)} Let $L$ be the minimal Quillen model of a simply connected complex of finite type $X$. Then, the cdga $\mathscr C^*(L)$ of cochains on $L$ has the homotopy type of  the Sullivan minimal model of $X$.
\end{corollary}
\begin{proof} It is enough (and in fact equivalent)  to check that the Sullivan realization of $\mathscr C^*(L)$ has the homotopy type of that of the Sullivan minimal model of $X$ which is its rationalization $X_\bq$. Indeed,
$$
\langle \mathscr C^*(L)\rangle_S \simeq \langle L\rangle\simeq\langle L\rangle_Q\simeq X_\bq.
$$
For the first equivalence one can use different arguments:  a direct proof is in  \cite[Theorem 8.1]{bfmt1}. Also one can consider the simplicial isomorphism $\langle \mathscr C^*(L)\rangle_S\cong\mc_\bullet(L)$, which is classical and easy to prove (see the Bibliographical Notes of \cite[Chapter 11]{bfmt0}), and then take into account that $\langle L\rangle$ is a deformation retract of $\mc_\bullet(L)$. The second equivalence is Theorem \ref{intro}.
\end{proof}

\section{Extending the Baues-Lemaire Conjecture}
Recall from \cite[\S3.2]{bfmt0} that the {\em minimal model} of a connected cdgl $M$ is a connected cdgl of the form $(\libc(V),d)$, in which $d$ is decomposable, together with a quasi-isomorphism
$$
(\libc(V),d)\stackrel{\simeq}{\longrightarrow} M.
$$
The {\em minimal model} of a connected simplicial set  $X$  is the minimal model of  $\lasu_X^a$ where $a$ is any of the $0$-simplices of $X$ regarded as a Maurer-Cartan element of $\lasu_X$, see \cite[\S8.4]{bfmt0}.

In this context, the Baues-Lemaire conjecture has already been extended to the connected case: if $L$ is the minimal Lie model of a connected simplicial set of finite type $X$ then, Theorem 10.8 of \cite{bfmt0} guarantees  that  $\varinjlim_n\quic^*(L/L^n)$ is a Sullivan model of $X$.

On the othe hand, the generalization of rational homotopy theory to non nilpotent spaces given by Pridham in \cite{prid}, by means of the pro-category of nilpotent, finite type, differential graded Lie algebras, also admits a deep extension of the Baues-Lemaire conjecture \cite[Corollary 4.41
and Remark 4.42]{prid}

Another way to interpret this conjecture in this extended scenario is via the category $\catedgl$ of {\em complete enriched} dgl's (edgl henceforth), see \cite{feha} for details. A connected dgl $L$ is complete enriched if there is a family $\{I_\alpha\}_{\alpha\in\cal J}$, where $\cal J$ is well ordered and directed, satisfying:

\smallskip
(i) $\alpha\ge \beta$ if and only if $I_\beta\subset I_\alpha$.

\smallskip
(ii) $L/I_\alpha$ is a finite type nilpotent dgl for each $\alpha\in\cal J$.

\smallskip
(iii) $L=\varprojlim_{\alpha\in\cal J} L/I_\alpha$.

\smallskip
Given an edgl $L$ one can define the {\em enriched cochains} by
$$
\quic_e^*(L)=\varinjlim_{\alpha\in\cal J} \quic^*(L/I_\alpha).
$$
By Theorem 4 of  \cite[\S15]{feha} the enriched cochains functors establishes an equivalence between the homotopy categories of edgl's and connected cdga's. In particular, the Sullivan model of any connected simplicial set can be regarded as the enriched cochains of some edgl.

We remark that the categories $\catedgl$ and $\catcdgl$ are different in general but given $L$ a dgl such that $L/[L,L]$ is of finite type, then $L$ is a cdgl if and only if it is an edgl. In particular $\lasu_\bullet$ is a cosimplicial edgl and for each edgl $L$ we may define its realization as the simplicial set
\begin{equation}\label{enr}
\langle L\rangle_e=\Hom_{\catedgl}(\lasu_\bullet,L).
\end{equation}
By the same observation, if $L$ is a cdgl such that $L/[L,L]$ is of finite type then $\langle L\rangle=\langle L\rangle_e$. However, in general, $\langle \,\cdot\,\rangle_e$
coincides with the Sullivan realization of the enriched cochains:
\begin{proposition} Let $L$ be any edgl. Then, $\langle \quic_e^*(L)\rangle_S\simeq\langle L\rangle_e$.
\end{proposition}
\begin{proof}
Since the Sullivan realization functor is a contravariant right adjoint $\langle \quic_e^*(L)\rangle_S\simeq\varprojlim_{\alpha\in\cal J}\langle\quic^*(L/I_\alpha)\rangle_S$. But, as each $L/I\alpha$ is of finite type, the above observation tells us that $\langle\quic^*(L/I_\alpha)\rangle_S\simeq\langle L/I_\alpha\rangle\simeq\langle L/I_\alpha\rangle_e$. Finally, by (\ref{enr}), $\langle\,\cdot\,\rangle_e$ commutes with limits and we conclude that
$\langle \quic_e^*(L)\rangle_S\simeq\langle L\rangle_e$.
\end{proof}

\begin{example} In the non finite type setting, the realization of a Sullivan model of a space which, as remarked above, is the enriched cochain functor of some edgl, may not coincide with the rationalization of the space. Or equivalently, with  the realization of the (unenriched) Quillen model of the space. For instance, let $X=\vee_{i\ge 1}S_i^2$ be a countably infinite wedge of $2$-spheres.  Its Quillen model is the simply connected (and hence complete) dgl $L=(\lib(V),0)$ where $V$ is a countably infinite dimensional vector spaces concentrated in degree $1$. Its realization $\langle L\rangle\simeq\langle L\rangle_Q$ is then homotopy equivalent to $\vee_{i\ge 1}(S_i^2)_\bq$.

On the other hand, let $\{v_i\}_{i\ge 1}$  be a basis of $V$, and denote $V(n)\subset V$ the vector space generated by the $v_i$, $i\geq n$. Then, the family of ideals  $\{I_n\}_{n\ge 1}$ of $L$, where $I_n = \cup_{j=0}^{n-1} \mathbb L^{\geq n-j}\bigl(V(j+1)\bigr)$, defines an edgl $L'=\varprojlim_n L/I_n$ whose realization $\langle L'\rangle_e$, as remarked above, is weakly homotopy equivalent to the realization of the Sullivan minimal model of $X$. That is,
$\varprojlim_n  \vee_{i\leq n} (S^2_i)_{\mathbb Q}$.

\end{example}

\bigskip
\noindent {\sc Institut de Math\'ematiques et Physique, Universit\'e Catholique de Louvain, Chemin du Cyclotron 2,
1348 Louvain-la-Neuve,
         Belgique}.

\noindent\texttt{yves.felix@uclouvain.be}

\medskip

\noindent{\sc Departamento de \'Algebra, Geometr\'{\i}a y Topolog\'{\i}a, Universidad de M\'alaga, 29080 M\'alaga, Spain.}

\noindent
\texttt{m\_fuentes@uma.es, aniceto@uma.es}
\end{document}